\documentclass[10pt]{amsart}

\usepackage[utf8]{inputenc}
\usepackage[english]{babel}
\usepackage{amssymb, amsmath, amsfonts}
\usepackage{xcolor}
\usepackage{hyperref}

\newtheorem{theo}{Theorem}[section]
\newtheorem{lem}[theo]{Lemma}
\newtheorem{cor}[theo]{Corollary}
\theoremstyle{definition}
\newtheorem{defi}[theo]{Definition}
\newtheorem{rem}[theo]{Remark}

\numberwithin{equation}{section}

\newcommand{\N}{\mathbb{N}}

\title{On weak$^*$-extensible subspaces of Banach spaces}

\author{G. Mart\'{\i}nez-Cervantes}
\address{Dpto. de Matem\'{a}ticas, Facultad de Matem\'{a}ticas, Universidad de Murcia, 30100 Espinardo (Murcia), Spain}
\email{gonzalo.martinez2@um.es}

\author{J. Rodr\'{i}guez}
\address{Dpto. de Ingenier\'{i}a y Tecnolog\'{i}a de Computadores,
Facultad de Inform\'{a}tica, Universidad de Murcia, 30100 Espinardo (Murcia), Spain}
\email{joserr@um.es}

\subjclass[2020]{Primary: 46B20, 46B50}

\keywords{Grothendieck space; $w^*$-extensible subspace; $G_\delta$-point; sequential compactness}

\begin{document}

\begin{abstract}
Let $X$ be a Banach space and $Y \subseteq X$ be a closed subspace. We prove that if the quotient $X/Y$ is weakly Lindelöf determined
or weak Asplund, then for every $w^*$-convergent sequence $(y_n^*)_{n\in \mathbb N}$ in~$Y^*$ there exist 
a subsequence $(y_{n_k}^*)_{k\in \mathbb N}$ and a $w^*$-convergent sequence $(x_k^*)_{k\in \mathbb N}$ in~$X^*$
such that $x_k^*|_Y=y_{n_k}^*$ for all $k\in \mathbb N$. As an application we obtain that $Y$ is Grothendieck 
whenever $X$ is Grothendieck and $X/Y$ is reflexive, which answers a question raised by Gonz\'{a}lez and Kania.
\end{abstract}

\maketitle

\section{Introduction}

Throughout this paper $X$ is a Banach space. We denote by $w^*$ the weak$^*$ topology on its (topological) dual~$X^*$.
The space $X$ is said to be {\em Grothendieck} if every $w^*$-convergent sequence in~$X^*$ is weakly convergent.
This property has been widely studied over the years, we refer the reader to the recent survey \cite{GonzalezKania}
for complete information on it. By a ``subspace'' of a Banach space we mean a closed linear subspace.
If $Y \subseteq X$ is a subspace, then: (i)~the quotient $X/Y$ is Grothendieck whenever~$X$ is Grothendieck, and (ii)~$X$
is Grothendieck whenever $Y$ and $X/Y$ are Grothendieck (see, e.g., \cite[2.4.e]{cas-gon}). In general, the property of being Grothendieck is not 
inherited by subspaces (for instance, $c_0$ is not Grothendieck while~$\ell_\infty$ is). However, this is the case for complemented subspaces or, more generally,
subspaces satisfying the following property:

\begin{defi}\label{defi:extensible}
A subspace $Y \subseteq X$ is said to be 
{\em $w^*$-extensible} in~$X$ if for every $w^*$-convergent sequence $(y_n^*)_{n\in \N}$ in~$Y^*$ there exist 
a subsequence $(y_{n_k}^*)_{k\in \N}$ and a $w^*$-convergent sequence $(x_k^*)_{k\in \N}$ in~$X^*$
such that $x_k^*|_Y=y_{n_k}^*$ for all $k\in \N$. 
\end{defi}

Indeed, it is easy to show that a Banach space is Grothendieck if (and only if) every $w^*$-convergent sequence in its dual admits 
a weakly convergent subsequence. Thus, a subspace $Y \subseteq X$ is Grothendieck whenever $X$ is Grothendieck
and $Y$ is $w^*$-extensible in~$X$. The concept of $w^*$-extensible subspace was studied in~\cite{cas-gon-pap,mor-wan,wan-alt}
(there the definition was given by replacing ``$w^*$-convergent'' by ``$w^*$-null''; both definitions are easily seen to be equivalent).
Note that every subspace is $w^*$-extensible in~$X$ whenever $B_{X^*}$ (the closed unit ball of~$X^*$, that we just call the ``dual ball'' of~$X$) is 
$w^*$-sequentially compact (cf. \cite[2.4.f]{cas-gon}). However, this observation does not provide new results on the stability of the Grothendieck property
under subspaces, because the only Grothendieck spaces having $w^*$-sequentially compact dual ball are the reflexive ones.

In this note we focus on finding sufficient conditions for the $w^*$-extensibility of a subspace~$Y\subseteq X$
which depend only on the quotient $X/Y$. Our motivation stems from the question (raised in \cite[Problem~23]{GonzalezKania}) of whether $Y$ is Grothendieck whenever 
$X$ is Grothendieck and $X/Y$ is reflexive. 
It is known that the separable injectivity of~$c_0$ (Sobczyk's theorem)
implies that if $X/Y$ is separable, then $Y$ is $w^*$-extensible in~$X$ in a stronger sense, namely, 
the condition of Definition~\ref{defi:extensible} holds without passing to subsequences (see, e.g., \cite[Theorem~2.3 and Proposition~2.5]{avi-alt-4}). 
Our main result is the following (see below for unexplained notation):

\begin{theo}\label{theo:Main}
Let $Y \subseteq X$ be a subspace. Suppose that $X/Y$ satisfies one of the following conditions:
\begin{enumerate}
	\item[(i)] Every non-empty $w^*$-closed subset of~$B_{(X/Y)^*}$ has a $G_\delta$-point
	(in the relative $w^*$-topology).
	\item[(ii)] ${\rm dens}(X/Y)<\mathfrak{s}$.
\end{enumerate}
Then $Y$ is $w^*$-extensible in~$X$.
\end{theo}

Given a compact Hausdorff topological space~$K$, a point $t\in K$ is called a {\em $G_\delta$-point} (in~$K$) if there is a sequence 
of open subsets of~$K$ whose intersection is $\{t\}$. Corson compacta have $G_\delta$-points 
(see, e.g., \cite[Theorem~14.41]{fab-ultimo}), and the same holds for any non-empty $w^*$-closed subset in the dual ball of a weak Asplund space 
(see, e.g., the proof of \cite[Theorem~2.1.2]{FabianDifferentiability}). Thus, we get the following corollary covering the case
when $X/Y$ is reflexive:

\begin{cor}\label{cor:classes}
Let $Y \subseteq X$ be a subspace. If $X/Y$ is weakly Lindelöf determined or weak Asplund, then 
$Y$ is $w^*$-extensible in~$X$.
\end{cor}

As an application of the above we get an affirmative answer to \cite[Problem~23]{GonzalezKania}:

\begin{cor}\label{cor:Grothendieck}
Let $Y \subseteq X$ be a subspace. If $X$ is Grothendieck and $X/Y$ is reflexive, then $Y$ is Grothendieck.
\end{cor}

As to condition~(ii) in Theorem~\ref{theo:Main}, recall that the density character of a Banach space~$Z$, denoted by~${\rm dens}(Z)$, is the
smallest cardinality of a dense subset of~$Z$. For our purposes, we just mention that
the {\em splitting number}~$\mathfrak{s}$ is the minimum of all cardinals~$\kappa$ for which there is a compact Hausdorff topological space
of weight~$\kappa$ that is not sequentially compact. In general, $\omega_1 \leq \mathfrak{s} \leq \mathfrak{c}$. So, under CH, cardinality 
strictly less than~$\mathfrak{s}$ just means countable. However, 
in other models there are uncountable sets of cardinality strictly less than~$\mathfrak{s}$. 
We refer the reader to~\cite{dou2} for detailed information on~$\mathfrak{s}$ and other cardinal characteristics of the continuum.

Any of the conditions in Theorem~\ref{theo:Main} implies that $B_{(X/Y)^*}$ is $w^*$-sequentially compact 
(see \cite[Lemma~2.1.1]{FabianDifferentiability} and note that the weight of $(B_{(X/Y)^*},w^*)$ coincides with ${\rm dens}(X/Y)$), 
which is certainly not enough to guarantee that $Y$ is $w^*$-extensible in~$X$ (see Remark~\ref{rem:BourgainSchlumprecht}). The ideas that we use in this paper are similar to those used by Hagler and Sullivan~\cite{HaglerSullivan} to study sufficient conditions for the dual ball of a Banach space to be $w^*$-sequentially compact. By the way,
as another consequence of Theorem~\ref{theo:Main} we obtain a generalization of \cite[Theorem~1]{HaglerSullivan} (see Corollary~\ref{cor:seqcompactness}).

The proof of Theorem~\ref{theo:Main} and some further remarks are included in the next section. We follow
standard Banach space terminology as it can be found in \cite{FabianDifferentiability} and~\cite{fab-ultimo}.

\section{Proof of Theorem~\ref{theo:Main} and further remarks}\label{section:proofs}

By a ``compact space'' we mean a compact Hausdorff topological space. The {\em weight} of a compact space~$K$, denoted by~${\rm weight}(K)$,
is the smallest cardinality of a base of~$K$. The following notion was introduced in~\cite{MP18}:

\begin{defi}\label{defi:CDE}
Let $L$ be a compact space and $K \subseteq L$ be a closed set. We say that $L$ is a \textit{countable discrete extension} of~$K$ if 
$L \setminus K$ consists of countably many isolated points.
\end{defi}
 
Countable discrete extensions turn out to be a useful tool to study twisted sums of~$c_0$ and $C(K)$-spaces, see \cite{avi-mar-ple} and~\cite{MP18}. 
As it can be seen in the proof of Theorem~\ref{theo:Main}, countable discrete extensions appear in a natural way when dealing with sequential properties and twisted sums. 
Lemma \ref{LemmaCDE} below isolates two properties of compact spaces
which are stable under countable discrete extensions, both of them implying sequential compactness. Nevertheless, sequential compactness itself is not stable under countable discrete extensions (see Remark~\ref{rem:BourgainSchlumprecht}).

\begin{lem}
\label{LemmaCDE}
Let $L$ be a compact space which 
is a countable discrete extension of a closed set $K \subseteq L$. Then:
\begin{enumerate}
\item[(i)] if every non-empty closed subset of~$K$ has a $G_\delta$-point (in its relative topology), then the same property holds for~$L$;
\item[(ii)] ${\rm weight}(L)={\rm weight}(K)$ whenever $K$ is infinite.
\end{enumerate}
Therefore, if either every non-empty closed subset of~$K$ has a $G_\delta$-point (in its relative topology) or ${\rm weight}(K)<\mathfrak{s}$, 
then $L$ is sequentially compact.
\end{lem}
\begin{proof} (i) Let $M \subseteq L$ be a non-empty closed set. If $M \subseteq K$, then $M$ has a $G_\delta$-point (in the relative topology) by hypothesis.
On the other hand, if $M \cap (L \setminus K)\neq \emptyset$, then $M$ contains a point which is isolated in~$L$ 
and so a $G_\delta$-point (in~$M$).

(ii) Let $R: C(L) \to C(K)$ be the bounded linear operator defined by $R(f)=f|_K$ for all $f\in C(K)$. Then $R$ is surjective
and $C(K)$ is isomorphic to~$C(L)/\ker R$. 
The fact that $L$ is a countable discrete extension of~$K$ implies that $\ker R$ is finite-dimensional or isometrically isomorphic to~$c_0$. 
In any case, $\ker R$ is separable and so 
$$
	{\rm dens}(C(K)) = {\rm dens}(C(L)). 
$$
The conclusion follows from the equality 
${\rm dens}(C(S))={\rm weight}(S)$, which holds for any infinite compact space~$S$ 
(see, e.g., \cite[Proposition 7.6.5]{Semadeni} or \cite[Exercise~14.36]{fab-ultimo}).

The last statement of the lemma follows from \cite[Lemma 2.1.1]{FabianDifferentiability} and \cite[Theorem~6.1]{dou2}, respectively.
\end{proof}

We will use below the well-known fact that ${\rm dens}(Z)={\rm weight}(B_{Z^*},w^*)$ for any Banach space~$Z$.

\begin{proof}[Proof of Theorem~\ref{theo:Main}]
Let $(y_n^*)_{n\in \N}$ be a $w^*$-convergent sequence in~$Y^*$. 
Without loss of generality, we can assume that $(y_n^*)_{n\in \N}$ is $w^*$-null and contained in~$B_{Y^*}$. Clearly, there is nothing 
to prove if $y_n^*=0$ for infinitely many $n\in \N$. 
So, we can assume further that $y_n^*\neq 0$ for all $n\in \N$.  
By the Hahn-Banach theorem, for each $n\in \N$ there is $z_n^* \in X^*$ with $z_n^*|_Y=y_n^*$ and $\|z_n^*\|=\|y_n^*\| \leq 1$.

Let $q:X\rightarrow X/Y$ be the quotient operator. It is well-known that its adjoint $q^*:(X/Y)^* \rightarrow X^*$ is an isometric isomorphism 
from $(X/Y)^*$ onto~$Y^\perp$. In addition, $q^*$ is $w^*$-to-$w^*$-continuous, hence $K:=B_{X^*}\cap Y^\perp=q^*(B_{(X/Y)^*})$ is $w^*$-compact.
Observe that for every $x^* \in X^* \setminus Y^\perp$ there exist $x\in Y$, $\alpha>0$ and $m\in \N$ 
such that $x^*(x) > \alpha > z_n^*(x)$ for every $n \geq m$. Hence $L:=K \cup \{z_n^*: n \in \N\} \subseteq B_{X^*}$
is $w^*$-closed (so that $L$ is $w^*$-compact) and each $z_n^*$ is $w^*$-isolated in~$L$ (bear in mind that $z_n^*|_Y=y_n^*\neq 0$).
Then $(L,w^*)$ is a countable discrete extension of $(K,w^*)$, with 
$(K,w^*)$ and $(B_{(X/Y)^*},w^*)$ being homeomorphic.
Bearing in mind that ${\rm dens}(X/Y)$ coincides with the weight of $(B_{(X/Y)^*},w^*)$, from Lemma~\ref{LemmaCDE} 
it follows that $L$ is sequentially compact and, therefore, $(z_n^*)_{n\in \N}$ admits a $w^*$-convergent subsequence. The proof is finished.
\end{proof}

If $X$ is Grothendieck and $Y \subseteq X$ is a subspace such that $X/Y$ is separable, then $Y$ is Grothendieck
(see \cite[Proposition~3.1]{gon-alt}). This fact can also be seen as a consequence of Corollary~\ref{cor:Grothendieck}, because
every separable Grothendieck space is reflexive. 

As an immediate application of Theorem~\ref{theo:Main} we also obtain an affirmative answer to
\cite[Problem~22]{GonzalezKania} (bear in mind that $\mathfrak{p}\leq \mathfrak{s}$, see, e.g., \cite[Theorem~3.1]{dou2}):

\begin{cor}\label{cor:Grothendieck-s}
Let $Y \subseteq X$ be a subspace such that ${\rm dens}(X/Y)<\mathfrak{s}$. If $X$ is Grothendieck, then $Y$ is Grothendieck.
\end{cor}

\begin{rem}\label{rem:alternate}
In fact, the previous corollary is a particular case of Corollary~\ref{cor:Grothendieck}. Indeed, on the one hand, 
the assumption that ${\rm dens}(X/Y)<\mathfrak{s}$ implies that $B_{(X/Y)^*}$ is $w^*$-sequentially compact.
On the other hand, the Grothendieck property is preserved by quotients
and any Grothendieck space with $w^*$-sequentially compact dual ball is reflexive (cf. \cite[Proposition~6.18]{avi-alt-4}). 
\end{rem}

\begin{rem}\label{rem:BourgainSchlumprecht}
In general, the $w^*$-sequential compactness
of~$B_{(X/Y)^*}$ is not enough to guarantee that a subspace $Y \subseteq X$ is $w^*$-extensible in~$X$. Indeed,
it is easy to check that {\em if both $B_{Y^*}$ and $B_{(X/Y)^*}$ are $w^*$-sequentially compact and $Y$
is $w^*$-extensible in~$X$, then $B_{X^*}$ is $w^*$-sequentially compact as well}
(see the proof of \cite[Proposition~6]{cas-gon-pap}).
On the other hand, there exists a Banach space $X$ such that $B_{X^*}$ is not $w^*$-sequentially compact
although $B_{(X/Y)^*}$ is $w^*$-sequentially compact for some separable subspace $Y \subseteq X$ (see \cite{HaglerSullivan}, cf. \cite[Section~4.8]{cas-gon}).
\end{rem}

Hagler and Sullivan proved in \cite[Theorem~1]{HaglerSullivan} that $B_{X^*}$ is $w^*$-sequentially compact whenever there is 
a subspace $Y \subseteq X$ such that $B_{Y^*}$ is $w^*$-sequentially compact and $X/Y$ has an equivalent G\^{a}teaux smooth norm.
Since every Banach space admitting an equivalent G\^{a}teaux smooth norm is weak Asplund (see, e.g., \cite[Corollary~4.2.5]{FabianDifferentiability}), 
the following corollary generalizes that result:

\begin{cor}
\label{cor:seqcompactness}
Let $Y \subseteq X$ be a subspace such that $B_{Y^*}$ is $w^*$-sequentially compact. 
If $X/Y$ satisfies any of the conditions in Theorem~\ref{theo:Main}, 
then $B_{X^*}$ is $w^*$-sequentially compact.
\end{cor}

\subsection*{Acknowledgements} The research is partially supported by {\em Agencia Estatal de Investigaci\'{o}n} [MTM2017-86182-P, grant cofunded by ERDF, EU] 
and {\em Fundaci\'on S\'eneca} [20797/PI/18]. The research of G. Mart\'inez-Cervantes was also partially supported 
by the European Social Fund and the Youth European Initiative under {\em Fundaci\'{o}n S\'{e}neca} [21319/PDGI/19].

\bibliographystyle{amsplain}

\end{document}